\documentclass{amsart}

\usepackage{url,graphicx}
\usepackage[dvips]{epsfig}
\usepackage{latexsym,color}
\usepackage{amscd,amssymb,verbatim}
\newcommand{\aut}{\text{\rm Aut}\,}

\newcommand{\bd}{{\textrm{Bd}\,}}
\newcommand{\be}{\begin{enumerate}}

\newcommand{\ca}{{\mathcal A}}

\newcommand{\chom}{\text{\tt Hom}_{\,0}}

\newcommand{\cn}{{\mathcal N}}

\newcommand{\conn}{\text{ conn}\,}
\newcommand{\cp}{{\mathcal P}}

\newcommand{\cs}{{\mathcal S}}

\newcommand{\da}{\Delta}

\newcommand{\dz}{{\mathbb Z}}

\newcommand{\ee}{\end{enumerate}}

\newcommand{\hra}{\hookrightarrow}

\newcommand{\id}{\text{id}}

\newcommand{\lk}{\text{lk}\,}

\newcommand{\lra}{\longrightarrow}

\newcommand{\mfq}{{\mathfrak Q}}

\newcommand{\nin}{\noindent}

\newcommand{\ovr}[1]{\overline{#1}}

\newcommand{\pr}{\noindent{\bf Proof. }}

\newcommand{\ra}{\rightarrow}

\newcommand{\sm}{\setminus}

\newcommand{\supp}{\text{\rm supp}\,}

\newcommand{\thom}{\text{\tt Hom}\,}

\newcommand{\ti}{\tilde }
\newcommand{\tn}{\text{\tt N}}

\newcommand{\wti}{\widetilde }

\newcommand{\zz}{{{\mathbb Z}_2}}

\newtheorem{thm}{Theorem}[section]
\newtheorem{df}  [thm]{Definition}

\newtheorem{crl} [thm]{Corollary}
\newtheorem{prop}[thm]{Proposition}
\newtheorem{conj}[thm]{Conjecture}
\newtheorem{rem} [thm]{Remark}

\numberwithin{equation}{section}

\numberwithin{figure}{section}

\title {Complexes of graph homomorphisms}

\author{Eric Babson and Dmitry N. Kozlov}

\date{\noindent
\today\\[0.05cm]
     \hskip15pt The second author acknowledges support by the
University of Washington, Seattle,
the Swiss National Science Foundation Grant PP002-102738/1, 
the University of Bern, and the Royal Institute of Technology, Stockholm.\\[0.05cm]
            }

\address{Department of Mathematics, University of Washington, Seattle,
  U.S.A.} \email{babson@math.washington.edu}

\address{current address: Department of Computer Science, Eidgen\"ossische Technische
Hochschule, Z\"urich, Switzerland}

\address{Department of Mathematics, Royal Institute of Technology,
Stockholm, Sweden.} \email{dkozlov@inf.ethz.ch}

\begin{document}

\begin{abstract}
  $\thom(G,H)$ is a~polyhedral complex defined for any two undirected
  graphs $G$ and~$H$. This construction was introduced by Lov\'asz to
  give lower bounds for chromatic numbers of graphs. In this paper we
  initiate the study of the topological properties of this class of
  complexes.
  
  We prove that $\thom(K_m,K_n)$ is homotopy equivalent to a wedge of
  $(n-m)$-dimensional spheres, and provide an enumeration formula for
  the number of the spheres.  As a corollary we prove that if for some
  graph $G$, and integers $m\geq 2$ and $k\geq -1$, we have
  $\varpi_1^k(\thom(K_m,G))\neq 0$, then $\chi(G)\geq k+m$; here
  $\zz$-action is induced by the swapping of two vertices in $K_m$, and
  $\varpi_1$ is the first Stiefel-Whitney class corresponding to this
  action.
  
  Furthermore, we prove that a~fold in the first argument of
  $\thom(G,H)$ induces a~homotopy equivalence. It then follows that
  $\thom(F,K_n)$ is homotopy equivalent to a~direct product of
  $(n-2)$-dimensional spheres, while $\thom(\ovr{F},K_n)$ is homotopy
  equivalent to a~wedge of spheres, where $F$ is an arbitrary forest
  and $\ovr{F}$ is its complement.
\end{abstract}

\maketitle


\section{Introduction}\label{sect_intr}

\subsection{Definition of the main object.} $\,$ \vspace{5pt}

For any graph $G$, we denote the set of its vertices by $V(G)$, and
the set of its edges by $E(G)$, $E(G)\subseteq V(G)\times V(G)$. All
the graphs in this paper are undirected, so $(x,y)\in E(G)$ implies
$(y,x)\in E(G)$. Unless otherwise specified, our graphs are finite and
may contain loops.

\begin{df}
  For two graphs $G$ and $H$, a {\bf graph homomorphism} from $G$ to
  $H$ is a~map $\phi:V(G)\rightarrow V(H)$, such that if $x,y\in V(G)$
  are connected by an edge, then $\phi(x)$ and $\phi(y)$ are also
  connected by an edge.
\end{df}

We denote the set of all homomorphisms from $G$ to $H$ by
$\chom(G,H)$.

\begin{df} \label{dfhom}
  $\thom(G,H)$ is a polyhedral complex whose cells are indexed by all
  functions $\eta:V(G)\rightarrow 2^{V(H)}\setminus\{\emptyset\}$,
  such that if $(x,y)\in E(G)$, then $\eta(x)\times\eta(y)\subseteq
  E(H)$.

  The closure of a~cell $\eta$ consists of all cells indexed by
  $\ti\eta:V(G)\rightarrow 2^{V(H)}\setminus\{\emptyset\}$, which
  satisfy $\ti\eta(v)\subseteq\eta(v)$, for all $v\in V(G)$.
\end{df}
\nin The set of vertices of $\thom(G,H)$ is precisely $\chom(G,H)$.
Since all cells of $\thom(G,H)$ are products of simplices, the
geometric realization of $\thom(G,H)$ is defined in a~straightforward
fashion.

On the intuitive level, one can think of each $\eta:V(G)\rightarrow
2^{V(H)}\setminus\{\emptyset\}$, satisfying the conditions of the
Definition~\ref{dfhom}, as associating non-empty lists of vertices of
$H$ to vertices of $G$ with the condition on this collection of lists
being that any choice of one vertex from each list will yield a~graph
homomorphism from $G$ to~$H$.

A direct geometric construction of $\thom(G,H)$ is as follows.
Consider the partially ordered set $P_{G,H}$ of all $\eta$ as in
Definition~\ref{dfhom}, with the partial order defined by
$\ti\eta\leq\eta$ if and only if $\ti\eta(v)\subseteq\eta(v)$, for all $v\in
V(G)$. Then the order complex $\Delta(P_{G,H})$ is a~barycentric
subdivision of $\thom(G,H)$. A~cell $\tau$ of $\thom(G,H)$ corresponds
to the union of all the simplices of $\Delta(P_{G,H})$ labeled by the
chains with the maximal element~$\tau$.

In this paper we study properties of the complexes $\thom(G,H)$.  More
specifically we compute the homotopy type of $\thom(G,H)$ for several
families of $G$ and $H$ and also derive some information about natural
finite group actions on these complexes.

\subsection{Historic motivation.} $\,$ \vspace{5pt}

A particularly frequently studied special case of a graph
homomorphism is that of a vertex coloring: for a graph $G$
a~vertex coloring of $G$ with $n$ colors is simply a~graph
homomorphism from $G$ to $K_n$. Here $K_n$ denotes an unlooped
complete graph on $n$ vertices, that is $V(K_n)=[n]$,
$E(K_n)=\{(x,y)\,|\,x,y\in [n], x\neq y\}$.

Historically, one was especially interested in the question of
existence of vertex colorings with a specified number of colors. From
this point of view, the minimal possible number of colors in a~vertex
coloring is of special importance. It is called the {\it chromatic
  number} of the graph, and is denoted by $\chi(G)$.

The Kneser conjecture was posed in 1955, see \cite{Knes}, and
concerned chromatic numbers of a specific family of graphs, later
called {\it Kneser graphs}. For $n,k\in\dz$, $n\geq 2$, $1\leq k\leq
n/2$, the Kneser graph $\Gamma_{k,n}$ is the graph whose vertices are
all $k$-subsets of $[n]$, and edges are all pairs of disjoint
$k$-subsets; here $1\leq k\leq n/2$.

In 1978 L.\ Lov\'asz solved the Kneser conjecture by finding geometric
obstructions of Borsuk-Ulam type to the existence of graph colorings.

\begin{thm}
  {\rm (Kneser-Lov\'asz, \cite{Knes,Lo}).} $\chi(\Gamma_{k,n})=n-2k+2$.
\end{thm}

To show the inequality $\chi(\Gamma_{k,n})\geq n-2k+2$ Lov\'asz
associated a~simplicial complex $\cn(G)$, called the
{\it~neighborhood complex}, to an~arbitrary graph $G$, and then
used the connectivity information of the topological space
$\cn(G)$ to find obstructions to the colorability of~$G$.

\begin{thm} \label{lothm}
  {\rm (Lov\'asz, \cite{Lo}).}  {\it Let $G$ be a graph, such that
    $\cn(G)$ is $k$-connected for some $k\in\dz$, $k\geq -1$, then
    $\chi(G)\geq k+3$.}
\end{thm}

The main topological tool which Lov\'asz employed was the Borsuk-Ulam
theorem. See also \cite{AFL} for the extension to hypergraphs, 
which used the generalization of the Borsuk-Ulam theorem from \cite{BSS}.

We shall define the complex $\cn(G)$ in Section~\ref{s4}, where we
shall also see that for any graph $G$ the complex $\cn(G)$ is
homotopy equivalent to $\thom(K_2,G)$. This fact leads one to
consider the family of $\thom$ complexes as a natural context in
which to look for further obstructions to the existence of graph
homomorphisms. 

Accordingly, Lov\'asz has made the following
conjecture, \cite{Lo2}. Let $C_m$ be a cycle with $m$ vertices:
$V(C_m)=\dz_m$, $E(C_m)=\{(x,x+1),(x+1,x)\,|\,x\in\dz_m\}$.

\begin{conj} \label{loconj}
  {\rm (Lov\'asz).}  {\it Let $G$ be a graph, such that
    $\thom(C_{2r+1},G)$ is $k$-connected for some $r,k\in\dz$, $r\geq
    1$, $k\geq -1$, then $\chi(G)\geq k+4$.}
\end{conj}

Our proof of Conjecture~\ref{loconj} was announced in~\cite{BK1}. 
The full version of the proof consists of 2 parts: the study of 
the important general properties of $\thom$ complexes, appearing
in this paper, and more detailed and specific spectral sequence 
calculations, appearing in~\cite{BK2}.
However, the general study undertaken in this paper contains more 
than the results which we later use for our spectral sequence
computations. 

There was a~more general conjecture, also due to Lov\'asz.
 \begin{conj} \label{loconj2}
  {\rm (Lov\'asz).}  {\it Let $T$, $G$ be two graphs, then
   \begin{equation}\label{eq:lo1}
    \chi(G)\geq\chi(T)+\conn\thom(T,G)+1.
   \end{equation}
   }
\end{conj}
Here, $\conn X$ is the connectivity of the topological space $X$, i.e., 
the maximum $k$ such that $X$ is $k$-connected. Note that 
Conjecture~\ref{loconj} is the special case of 
the Conjecture~\ref{loconj2} corresponding to $T=C_{2r+1}$.

Hoory and Linial, \cite{HL}, gave a~counterexample to 
the Conjecture~\ref{loconj2}. In their counterexample $G=K_5$, 
and $T$ is a~graph with 9 vertices and 22 edges. Furthermore, 
$\chi(T)=5$, but $\thom(T,G)$ is connected, showing that 
the equation \eqref{eq:lo1} is false in general.

In this paper we show that Conjecture~\ref{loconj2} is true for $T=K_m$ 
(Lov\'asz himself proved \eqref{eq:lo1} for $T=K_2$).
More specifically, let $\zz$ act on $K_m$ for $m\geq 2$, by swapping 
the vertices 1 and 2 and fixing the vertices $3,\dots,m$. Since 
the graph homomorphism flips an~edge, it induces a~free $\zz$-action 
on $\thom(K_m,G)$, for an arbitrary graph $G$ without loops. 

For a CW complex $X$ on which $\zz$ acts freely let $\varpi_1(X)$
denote its {\it first Stiefel-Whitney class}, see \cite{tD}. As
a~corollary of our computations, we prove the following theorem.
\begin{thm}\label{cothm}
  Let $G$ be a graph, and let $m,k\in\dz$, such that $m\geq 2$, $k\geq
  -1$. If $\varpi_1^k(\thom(K_m,G))\neq 0$, then $\chi(G)\geq k+m$.
\end{thm}
Note that if a $\zz$-space $X$ is $k$-connected, then there exists
a~$\zz$-map $S^{k+1}_a\ra X$. The functoriality of Stiefel-Whitney
classes and the fact that $\varpi_1^{k+1}(S^{k+1}_a)\neq 0$ imply
that $\varpi_1^{k+1}(X)\neq 0$. Therefore Theorem~\ref{cothm} implies
Conjecture~\ref{loconj2} for $T=K_m$.

\subsection{Plan of the paper.} $\,$
\vspace{5pt}

In Section 2 we define notations, describe the category of graphs and
graph homomorphisms, and give several examples of $\thom$ complexes.
Furthermore, we list many simple, but fundamental properties of the
$\thom$ construction.

In Section 3 we describe two results from topological combinatorics
which we need for our arguments: a~proposition from Discrete Morse
theory, and a~Quillen-type result.

In Section 4 we see first that in general $\thom(K_2,G)$ is homotopy
equivalent to the neighborhood complex $\cn(G)$, implying in
particular that $\thom(K_2,K_n)\simeq S^{n-2}$. We observe that in
fact $\thom(K_2,K_n)$ is a~boundary complex of a polytope, on which
the natural $\zz$-action on the first argument, induces an antipodal
action. In the subsection~\ref{ss4.3} we prove the central result of
this section, namely we show that $\thom(K_m,K_n)$ is homotopy 
equivalent to a~wedge of $(n-m)$-dimensional spheres, and provide 
an~enumeration formula for the number of the spheres. 
As a~corollary we derive the Theorem~\ref{cothm}.

In Section 5 we prove that a~fold in the first argument of $\thom(G,H)$
induces a~homotopy equivalence. As a~corollary, we show for 
an~arbitrary forest $F$, that $\thom(F,K_n)$ is homotopy equivalent 
to a~direct product of $(n-2)$-dimensional spheres, while 
$\thom(\ovr{F},K_n)$ is homotopy equivalent to a~wedge of spheres. 
For an arbitrary tree $T$ with $\zz$-action we describe the 
$\zz$-homotopy type of $\thom(T,K_n)$.

To conclude the introduction, we refer the reader to a~recent survey 
of the previous studies of other complexes related to graph colorings,
see~\cite{MZ}.

\vspace{5pt}

\nin {\bf Acknowledgments.} We would like to thank L\'aszl\'o Lov\'asz
for insightful discussions, Sonja \v{C}uki\'c for a~remark, and the
anonymous referee for useful suggestions which helped us 
to improve the presentation significantly.

\section{Basic facts about $\thom$ complexes.}

\subsection{Terminology.} $\,$
\vspace{5pt}

\nin $\circ$ For a~graph $G$ we distinguish between looped and
unlooped complements, namely we let ${\ovr G}^{\text{\,l}}$ be the graph
defined by
$$V({\ovr G}^{\text{\,l}})=V(G),\,\,E({\ovr G}^{\text{\,l}})=(V(G)\times
V(G))\sm E(G),$$ 
while $\ovr G$ is the graph defined by
$$V(\ovr G)=V(G),\,\,E(\ovr G)=\{(x,y)\in V(G)\times V(G)\,|\, x\neq
y, (x,y)\notin E(G)\}.$$

\nin $\circ$ For a~graph $G$ and $S\subseteq V(G)$ we denote by $G[S]$
the graph on the vertex set $S$ induced by $G$, that is $V(G[S])=S$,
$E(G[S])=(S\times S)\cap E(G)$. For $S\subseteq V(G)$ we set $G-S$ to
be the graph $G[V(G)\setminus S]$. For $v\in V(G)$ we shall sometimes
simply write $G-v$ instead of $G-\{v\}$.

\nin $\circ$ For a~graph $G$ and $A\subseteq V(G)$, let $\tn(A)=\{w\in
V(G)\,|\,(v,w)\in E(G),\,\forall v\in A\}$ denote the set of all
common neighbors of the vertices of $A$. In particular,
$\tn(\emptyset)=V(G)$, and $\tn(v):=\tn(\{v\})$ is simply the set of
all neighbors of $v$, with the convention being that $v$ is its own
neighbor if and only if $(v,v)\in E(G)$. If needed, we will also
specify the graph by writing $\tn_G(A)$.

\nin $\circ$ For two arbitrary graphs $G$ and $H$ we let $G\times H$
denote the {\it direct product} of $G$ and $H$:
\begin{multline*} V(G\times H)=V(G)\times V(H),\,\,
E(G\times H)=\{((x,y),(\ti x,\ti y))\,| \\
\,x,\ti x\in V(G),y,\ti y\in V(H),(x,\ti x)\in E(G),(y,\ti y)\in
E(H)\}.
\end{multline*}

\nin $\circ$ For two arbitrary graphs $G$ and $H$ we let $G\coprod H$
denote the {\it disjoint union} of $G$ and $H$.

\nin $\circ$ For $n\in\dz$, $n\geq 1$, we let $L_n$ denote the graph
defined by $V(L_n)=[n]$, $E(L_n)=\{(x,y)\,|\,|x-y|=1\}$.

\nin $\circ$ Let $\mfq$ be the graph defined by $V(\mfq)=[2]$,
$E(\mfq)=\{(1,2),(2,1),(1,1)\}$.

\nin $\circ$ For an arbitrary graph $G$, we let $G^o$ denote the
{\it~loop completion} of $G$, that is $V(G^o)=V(G)$,
$E(G^o)=E(G)\cup\{(v,v)\,|\,v\in V(G)\}$.

\nin $\circ$ For a~polyhedral complex $K$ we let $\cp(K)$ denote its
{\it face poset}, that is a~partially ordered set of the faces ordered
by inclusion.

\nin $\circ$ For any finite category $C$ (in particular a finite
poset) we denote by $\Delta(C)$ the realization of the nerve of that
category.

\nin $\circ$ For a~poset $P$ we let $\bd(P)$ denote the barycentric
subdivision of $P$, that is the poset of all the chains in the given
poset ordered by inclusion. For a~polyhedral complex $K$ we let
$\bd(K)$ denote the barycentric subdivision of $K$. Clearly,
$\bd(K)=\Delta(\cp(K))$, and $\cp(\Delta(P))=\bd(P)$.

\nin $\circ$ For any finite poset $P$, we let $P^{op}$ denote the
finite poset which has the same set of elements as $P$, but the
opposite partial order. Also, for any finite poset $P$, whenever
a~subset of the elements of $P$ is considered as a~poset, the partial
order is taken to be induced from~$P$.

\nin $\circ$ {\bf Top} is a~category having topological spaces as
objects, and continuous maps as morphisms.

\subsection{The category {\bf Graphs}.} $\,$
\vspace{5pt}

It is an easy check that a~composition of two graph homomorphisms is
again a~graph homomorphism. We denote a~composition of
$\phi\in\chom(G,H)$ and $\psi\in\chom(H,K)$ by
$\psi\circ\phi\in\chom(G,K)$.

Since the composition is associative and since for any graph $G$
we have a~unique identity homomorphism in $\chom(G,G)$ we can
define a~category {\bf~Graphs} as the one having graphs as
objects, and graph homomorphisms as morphisms.

One can check that the direct product of graphs is a~categorical
{\it~product} in {\bf Graphs}, while the disjoint union of graphs
is a~categorical {\it~coproduct} in {\bf Graphs}.

Note that with the above notations $K_1^{o}$ is a~graph consisting of
one vertex and one loop, it is the {\it terminal object} of
{\bf~Graphs}. The empty graph is the {\it initial object} of
{\bf~Graphs}.

\subsection{Examples of $\thom$ complexes.} $\,$ \label{ss2.3}
\vspace{5pt}

To start with, we have various trivial cases:

\nin $\circ$ $\thom(K_1,H)$ is a simplex with $|V(H)|$ vertices;

\begin{figure}[hbt]
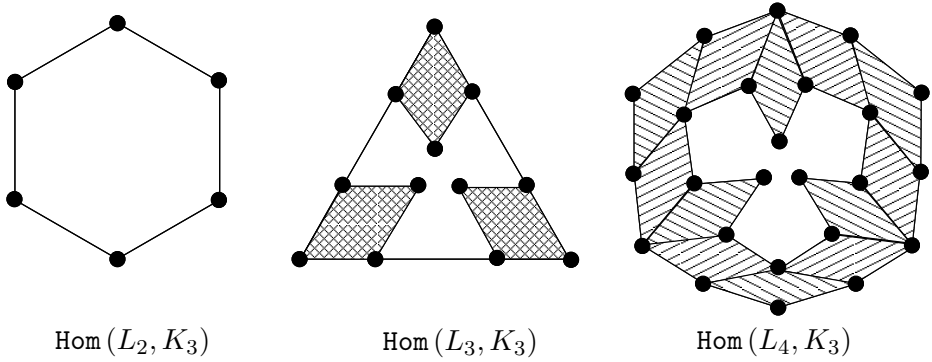

\begin{center}
  \begin{picture}(0,0)%
    \includegraphics{homlk.pstex}%
  \end{picture}%
  \input{homlk.pstex_t}%
  
\end{center}
\caption{Complexes of 3-colorings of strings.}
\label{fig:homlk}
\end{figure}

\nin $\circ$ $\thom(H,K_1)=\emptyset$, unless $E(H)=\emptyset$, in
which case $\thom(H,K_1)$ is a~point;

\nin $\circ$ more generally, $\thom(G,H)=\emptyset$ if
$\chi(G)>\chi(H)$;

\nin $\circ$ $\thom(K_1^o,H)$ is a simplex with vertices indexed by
the looped vertices of~$H$;

\nin $\circ$ $\thom(H,K_1^o)$ is a point, as mentioned above;

\nin $\circ$ $\thom(G,K_n^o)$ is a~direct product of $|V(G)|$
simplices, each simplex having $n$ vertices;

\nin $\circ$ $\thom(G,K_2)=\emptyset$ if $G$ is not bipartite; it
consists of $2^c$ points, if $G$ bipartite and has $c$ connected
components;

\nin $\circ$ $\thom(C_{2r+1},C_{2p+1})=\emptyset$ if and only if $r<p$; 

\nin $\circ$ $\thom(C_{2r+1},C_{2r+1})$ is a disjoint union of
$4r+2$ points, for $r\geq 1$;

\nin $\circ$ $\thom(C_{2r+1},C_{2r-1})$ is a disjoint union of
two cycles, each of length $4r^2-1$, for $r\geq 2$;

\begin{figure}[hbt]
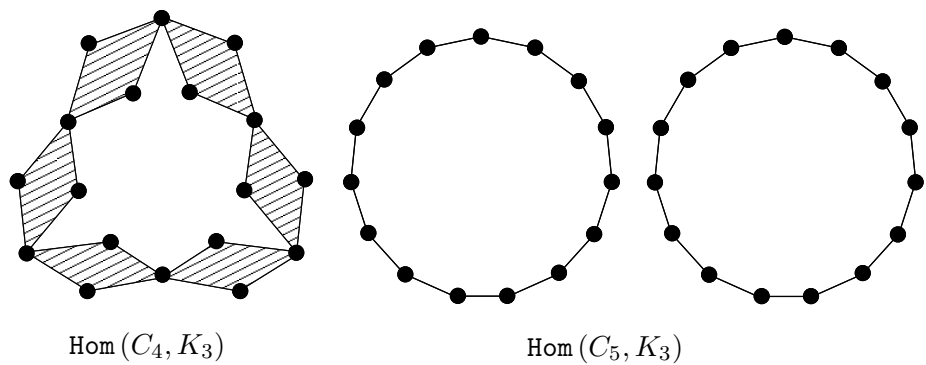

\begin{center}
  \begin{picture}(0,0)%
    \includegraphics{homc45k3.pstex}%
  \end{picture}%
  \input{homc45k3.pstex_t}%
  
\end{center}
\caption{Complexes of 3-colorings of 4- and 5-cycles.}
\label{fig:homc45k3}
\end{figure}

\begin{figure}[hbt]
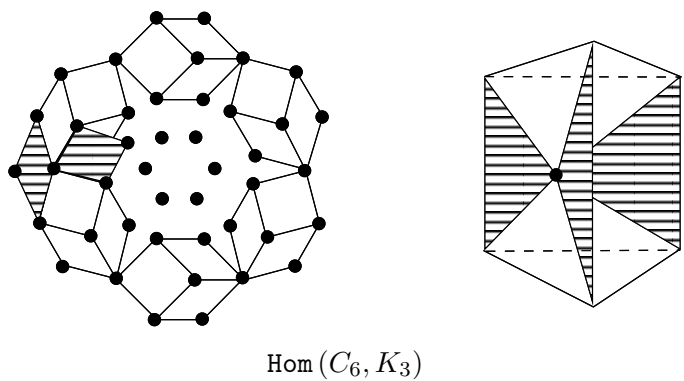

\begin{center}
  \begin{picture}(0,0)%
    \includegraphics{homc6k3.pstex}%
  \end{picture}%
  \input{homc6k3.pstex_t}%
  
\end{center}
\caption{Complexes of 3-colorings of a 6-cycle.}
\label{fig:homc6k3}
\end{figure}

\nin $\circ$ $\thom(C_6,K_3)$ consists of 6 isolated points, 6 solid
cubes and 18 squares connected as shown on the
Figure~\ref{fig:homc6k3}. The left part of Figure~\ref{fig:homc6k3} is
incomplete for the purpose of visualizing, it shows the 6 points, 6
cubes and some of the squares. The right part shows the link of each
of the 6 vertices, where two of the cubes touch. The closed star of
such a vertex consists of 2 solid cubes and 3 squares.

\nin $\circ$ $\thom(K_n,K_n)$ is a disjoint union of $n!$ points;

\nin $\circ$ $\thom(K_{n-1},K_n)$ is the Cayley graph of $\cs_n$
with the set of generators consisting of $n-1$ transpositions
$\{(a,n)\,|\,a=1,\dots,n-1\}$. Indeed, every vertex of
$\thom(K_{n-1},K_n)$ is an injection $\iota:[n-1]\ra[n]$, which
can be identified with a permutation of $[n]$ by writing out the
values of $\iota$ and then writing the missing element of $[n]$ in
the last position. An edge is a changing of one arbitrary value of
$\iota$, say $\iota(a)$, to the missing value, which is precisely
the same as acting with the transposition $(a,n)$ on the
corresponding permutation.

\begin{figure}[hbt]
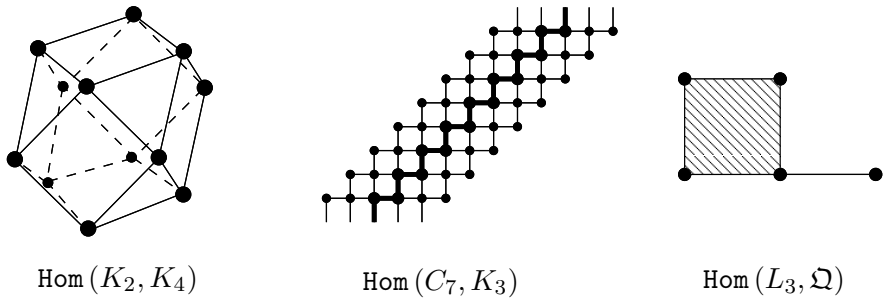

\begin{center}
  \begin{picture}(0,0)%
    \includegraphics{hom3.pstex}%
  \end{picture}%
  \input{hom3.pstex_t}%
  
\end{center}
\caption{Further examples of $\thom$-complexes.}
\label{fig:hom3}
\end{figure}

\nin $\circ$ $\thom(K_2,K_4)$ is the full 2-skeleton of the 3-cell
depicted on Figure~\ref{fig:hom3}.

\nin $\circ$ $\thom(C_7,K_3)$ is homeomorphic to a disjoint union of
two M\"obius bands. The local structure of each M\"obius band is shown
on the Figure~\ref{fig:hom3}. The middle cycle which is painted bold
has length 21 in each band, and all visible squares on the picture are
filled with 2-cells.

\nin $\circ$
It is not difficult to count the number of connected components
of $\thom(C_t,K_3)$. Denote this number $c_t$, the general formula is 
$$c_t=\begin{cases}
\lfloor(t+1)/3\rfloor,& \text{ if }3\hspace{-2.5pt}\not{|}\,\, t,\\
t/3 +5,& \text{ if }3\,\,{|}\,\, t,
\end{cases}$$
for $t\geq 3$. The crucial fact for deriving this formula for $c_t$ is
to notice that the connected components of $\thom(C_t,K_3)$ are
indexed with the number of times $C_t$ can wind around the triangle
$K_3$, with the sign encoding the direction.

\begin{rem}
As the examples above show, for the considered values of $m$ and $n$,
the spaces $\thom(C_m,C_n)$ are either empty or consist of several
connected components, with each component either being a point, or
homotopy equivalent to a~circle. Recently, this fact has been proved
for all values of $m$ and $n$, see~\cite{CK1}.
\end{rem}

\nin $\circ$ Let $V_{n,k}$ denote the Stiefel manifold
of all orthonormal $k$-frames in ${\mathbb R}^n$. 
Csorba, \cite{Cs}, has made the following conjecture:
\begin{conj}
(Csorba). $\thom(C_5,K_n)$ is homeomorphic to $V_{n-1,2}$, 
for all $n\geq 1$. 
\end{conj}
The cases $n=1,2$ are tautological, as both spaces are empty.
The example above verifies the case $n=3$: 
$\thom(C_5,K_3)\cong S^1\coprod S^1$. Several cases, including
$n=4$ have been recently verified by Csorba and Lutz,~\cite{CL}.

\nin $\circ$  Note that for an arbitrary $G$, $\thom(G,\mfq)$ can
be interpreted as a~cubical cone over the independence complex of
$G$; recall that the independence complex of $G$ is the simplicial
complex consisting of all independent sets of $G$. When saying
{\it cubical cone} we mean the following construction: given an
arbitrary simplicial complex $\da$, add an extra vertex $a$, and
for each simplex $\sigma\in\da$ with $d$ vertices span
a~$d$-dimensional cube $K_\sigma$ with $a$ being a~vertex of
$K_\sigma$ and $\sigma$ forming the link of $a$ in $K_\sigma$.

Note that $\thom(G,K_3)$ is cubical for any graph $G$ having no
isolated vertices. By a~theorem of Gromov, see~\cite{BH},
$\thom(G,K_3)$ allows metric with non\-positive curvature if and only
if the link of every vertex is a flag complex (which means that each
link is the clique complex of its 1-skeleton).

For any $\varphi\in\chom(G,H)$, we say that $\varphi$ has a~cubical
neighborhood if $\varphi$ does not belong to any simplex with more
than 2 vertices.

\begin{prop} \label{cubprop}
 If $\varphi\in\chom(G,H)$ has a~cubical neighborhood, then
$\lk_{\thom(G,H)}\varphi$ is a~flag complex.
\end{prop}

\pr Set $L=\lk_{\thom(G,H)}\varphi$, i.e., the simplicial complex
whose face poset is $\thom(G,H)_{\geq\varphi}$. For $v\in V(G)$, set
$$A_\varphi(v)=\tn(\bigcup_{w\in\tn(v)}\varphi(w)).$$
Since $\varphi$
has a~cubical neighborhood, we have $|A_\varphi(v)|\in\{1,2\}$, for
any $v\in V(G)$. Let $M(\varphi)\subseteq V(G)$ be the set of all
vertices $v$ with $|A_\varphi(v)|=2$.

Clearly $L$ has $M(\varphi)$ as the set of vertices. Furthermore,
$\sigma\subseteq M(\varphi)$, such that $|\sigma|\geq 2$, is a~simplex
in $L$, if and only if, for any two $a,b\in\sigma$, and any $x\in
A_\varphi(a)$, $y\in A_\varphi(b)$, we have $(x,y)\in E(H)$. Since
this is a~local condition depending only on the pair $(a,b)$, we
conclude that $L$ is a~flag complex.  \qed

It follows that the cubical complex $\thom(G,K_3)$ always allows
a~metric with non\-positive curvature. Moreover, for any
$\varphi\in\chom(G,K_3)$, the proof of the Proposition~\ref{cubprop}
yields that $\lk_{\thom(G,K_3)}\varphi$ is the independence complex
of~$G[M(\varphi)]$.

\subsection{General properties of $\thom$ complexes.} $\,$\label{ss2.4}
\vspace{5pt}

\nin (1) For any two graphs $G$ and $H$, $\thom(G,H)$ is a~regular CW
complex. \vspace{5pt}

\nin (2) Cells of $\thom(G,H)$ are direct products of simplices. More
specifically, each $\eta$ as in the Definition~\ref{dfhom} is
a~product of $|V(G)|$ simplices, having dimensions $|\eta(x)|-1$, for
$x\in V(G)$. Thus $\dim\eta=\sum_{x\in V(G)}|\eta(x)|-|V(G)|$.
\vspace{5pt}

\nin (3) For any three graphs $G$, $H$, and $K$, we have
\[\thom(G\coprod H,K)=\thom(G,K)\times\thom(H,K),\]
and, if $G$ is connected, and $G\neq K_1$, then also
\[\thom(G,H\coprod K)=\thom(G,H)\coprod\thom(G,K),\]
where the equality denotes isomorphism of polyhedral complexes.

The first formula is obvious. To see the second one, note that for
$\eta:V(G)\ra 2^{V(H)\cup V(K)}\sm\{\emptyset\}$, and $x,y\in V(G)$,
such that $(x,y)\in E(G)$, if $\eta(x)\cap V(H)\neq\emptyset$, then
$\eta(y)\subseteq V(H)$, which under assumptions on $G$ implies that
$\bigcup_{x\in V(G)}\eta(x)\subseteq V(H)$.
\vspace{5pt}

\nin (4) $\thom(H,-)$ is a covariant, while $\thom(-,H)$ is
a~contravariant functor from {\bf Graphs} to {\bf Top}.

If $\phi\in\chom(G,G')$, then we shall denote the cellular maps
induced by composition as $\phi^H:\thom(H,G)\ra\thom(H,G')$ and
$\phi_H:\thom(G',H)\ra\thom(G,H)$. \vspace{5pt}

\nin (5) The map induced by composition
\[\thom(G,H)\times\thom(H,K)\lra\thom(G,K)\]
is a~topological map.
\vspace{5pt}

\nin (6) Obviously, it is difficult to decide in general whether
$\thom(G,K_n)$ is non-empty, let alone $k$-connected. It is certainly
non-empty if the valency of each vertex is at most $n-1$. The
following fact is true in general.
\begin{prop} \label{pr:conn} 
  Let $G$ be any graph.  If the maximal valency of $G$ is equal to
  $d$, then $\thom(G,K_n)$ is connected, for all $n\geq d+2$. 
\end{prop}
\pr Assume $\thom(G,K_n)$ is not connected. Choose
$\phi,\psi\in\chom(G,K_n)$, such that $\psi$ and $\phi$ belong to
different connected components, and $\phi(v)=\psi(v)$ for the maximal
possible number of vertices. Pick $u$, such that $\phi(u)\neq\psi(u)$.
If $\psi(u)$ cannot be changed to $\phi(u)$, that is, if $\tau:V(G)\ra
V(H)$, defined by $\tau(x)=\psi(x)$ for $x\neq u$, $\tau(u)=\phi(u)$,
is not a~graph homomorphism, then there exists a~vertex $w$, such that
$(u,w)\in E(G)$, and $\psi(w)=\phi(u)\neq\phi(w)$. 

Since the valency of $w$ is at most $n-2$, we can change $\psi(w)$ to
something else, without changing the number of vertices on which
$\psi$ and $\phi$ coincide. Once this is done for each such neighbor
of $u$, we can change $\psi(u)$ to $\phi(u)$, thereby increasing the
number of vertices on which $\psi$ and $\phi$ coincide, hence
obtaining a contradiction to the choice of $\psi$ and $\phi$. \qed

This result motivates the following conjecture.

\begin{conj}\label{conj1}
Let $G$ be any graph.  If the maximal valency of $G$ is equal 
to~$d$, then $\thom(G,K_n)$ is $k$-connected, for all integers 
$k\geq -1$, $n\geq d+k+2$. 
\end{conj}

Proposition~\ref{pr:conn} corresponds to the case $k=0$; 
the case $k=-1$ is obviously true since the graph of maximal 
valency $d$ can be colored with $d+1$ colors.

\section{Tools from topological combinatorics}

\subsection{Discrete Morse theory.} $\,$
\vspace{5pt}

For a poset $P$ with the covering relation $\succ$, we define
a~{\it~partial matching} on $P$ to be a~set $S\subseteq P$, and
an~injective map $\mu:S\ra P\sm S$, such that $\mu(x)\succ x$, for all
$x\in S$. The elements of $P\sm(S\cup\mu(S))$ are called critical.

The next proposition is a~special case, which will be sufficient for
our purposes, of a~more general result proved by R.\ Forman,
see~\cite{Fo}.

\begin{prop}\label{DMT}
  Let $\da$ be a~regular CW complex and $\da'$ a~subcomplex of $\da$,
  then the following are equivalent:

  a) there is a~sequence of collapses leading from $\da$ to $\da'$;

  b) there is a~partial matching $\mu$ on $\cp(\da)$ with the set of
  critical cells being $\cp(\da')$, such that there is no sequence 
  $x_1,\dots,x_t\in\cp(\da)\sm\cp(\da')$, $t\geq 2$, such that
  $\mu(x_1)\succ x_2, \mu(x_2)\succ x_3,\dots,\mu(x_t)\succ x_1$ (such
  matching is called acyclic).
\end{prop}
\pr See \cite[Proposition 5.4]{Koz2}.
\qed\vspace{5pt}

Proposition~\ref{DMT} is a~part of the Discrete Morse theory;
\cite{BBL,Fo,Koz1,Koz2} are just some of the references where it has
been studied and used.

\subsection{A Quillen-type result.} $\,$
\vspace{5pt}

In this subsection we prove a~Quillen-type result which, given a~poset
map $\phi$ satisfying certain conditions, provides us with some
topological information about the induced simplicial map~$\da(\phi)$.

\begin{prop} \label{prop_ABC}
Let $\phi:P\ra Q$ be a map of finite posets. Consider a list of
possible conditions on $\phi$.

\nin {\rm Condition $(A)$.} For every $q\in Q$, $\da(\phi^{-1}(q))$ is
contractible.

\nin {\rm Condition $(B)$.} For every $p\in P$ and $q\in Q$ with
$p\in\phi^{-1}(Q_{\geq q})$ the poset $\phi^{-1}(q)\cap P_{\leq p}$
has a maximal element. In this case we denote this maximal element by
$\max(p,q)$.

\nin {\rm Condition $(B^{op})$.} Let $\phi^{op}:P^{op}\ra Q^{op}$ be
the poset map induced by $\phi$. We require that $\phi^{op}$ satisfies
Condition~$B$. In this case we denote the minimal element of
$\phi^{-1}(q)\cap P_{\geq p}$ by $\min(p,q)$.

Then
\begin{enumerate}
\item [(1)] If $\phi$ satisfies $(A)$ and either $(B)$ or $(B^{op})$,
  then $\Delta(\phi)$ is a~homotopy equivalence.

\item [(2)] If $\phi$ satisfies $(B)$ and $(B^{op})$, and $Q$ is
  connected, then for any $q,q'\in Q$ we have
  $\da(\phi^{-1}(q))\simeq\da(\phi^{-1}(q'))$.  Furthermore, we have
  a~fibration homotopy long exact sequence:
\begin{equation} \label{eq:piseq}
\dots\lra\pi_i(\da(\phi^{-1}(q)))\lra\pi_i(\da(P))\lra\pi_i(\da(Q))\lra\dots
\end{equation}
\end{enumerate}
\end{prop}
\pr Consider the poset map $\bd\phi:\bd P\ra\bd Q$, which maps
$\rho\in\bd P$, $\rho=(\alpha_1>\dots>\alpha_t)$ to
$\{\phi(\alpha_1),\dots,\phi(\alpha_t)\}$. Since $\phi$ is
order-preserving, the last set is totally ordered, and thus can be
interpreted as a~chain in~$Q$.

We set $\phi^{-1}(\gamma):=\bigcup_{i=1}^t\phi^{-1}(\alpha_i)$ and
view it as a~subposet of~P. Note that
\begin{equation} \label{eq:bdgam}
(\bd\phi)^{-1}(\bd Q_{\leq\gamma})=\bd(\phi^{-1}(\gamma)).
\end{equation}

First we show (1). Because of the symmetry, we restrict our
consideration to the case when $\phi$ satisfies conditions $(A)$ and
$(B^{op})$. By Quillen's theorem~A, see \cite[p.\ 85]{Qu}, it is
enough to show that $\da((\bd\phi)^{-1}(\bd Q_{\leq\gamma}))$ is
contractible for any $\gamma\in\bd Q$. By \eqref{eq:bdgam} it is
enough to show that $\da(\phi^{-1}(\gamma))$ is contractible for any
$\gamma\in\bd Q$. We use induction on the length of the chain
$\gamma=(\alpha_1>\dots>\alpha_t)$. When $t=1$, this is precisely
condition~$(A)$, so we assume that $t\geq 2$.

Define $\xi:\phi^{-1}(\gamma)\ra\phi^{-1}(\alpha_1)$, by
$\xi(p)=\min(p,\alpha_1)$, for $p\in\phi^{-1}(\gamma)$. This is
well-defined since $\phi(p)\leq\alpha_1$. Note that
\begin{enumerate}
\item[1)] $\xi^2=\xi$, since $\xi|_{\phi^{-1}(\alpha_1)}=\id$;
\item[2)] $\xi(p)\geq p$, by the definition of $\min(p,\alpha_1)$;
\item[3)] $\xi$ is order-preserving. Indeed, take
  $p,p'\in\phi^{-1}(\gamma)$, such that $p>p'$. Then, on one hand
  $\xi(p)\geq p>p'$, on the other hand $\phi(\xi(p))=\alpha_1$, hence,
  by the definition of $\min(p',\alpha_1)$, we have
  $\xi(p)\geq\xi(p')$.
\end{enumerate}
This means that $\xi$ is a~closure map, hence $\da(\xi)$ is homotopy
equivalence, see~\cite[Corollary 10.12]{Bj}. It follows by induction
that $\da(\phi^{-1}(\gamma))$ is contractible for any $\gamma\in\bd
Q$.

Next we prove (2). Let $\gamma,\ti\gamma\in\bd Q$, such that
$\gamma>\ti\gamma$. We want to show that the inclusion map
$i:\phi^{-1}(\ti\gamma)\hookrightarrow\phi^{-1}(\gamma)$ induces
a~homotopy equivalence of the order complexes. Set $\gamma'=\gamma\cap
Q_{\geq\min\ti\gamma}$. Then $\min\ti\gamma=\min\gamma'$,
$\max\gamma=\max\gamma'$, and $\gamma\geq\gamma'$.

Consider the sequence of inclusion maps $\phi^{-1}(\max\gamma)
\stackrel{i_1}\hookrightarrow\phi^{-1}(\gamma')
\stackrel{i_2}\hookrightarrow\phi^{-1}(\gamma)$, and let
$\xi:\phi^{-1}(\gamma)\rightarrow\phi^{-1}(\max\gamma)$ be the map
defined above. By the argument for the part (1) we know that pairs
$(i_1,\xi)$ and $(i_2\circ i_1,\xi)$ induce homotopy equivalences
of the order complexes. It follows that the pair $(i_2,\xi)$ also
induces a~homotopy equivalence, since
$$\da(i_2)\circ\da(\xi)=\da(i_2)\circ\da(i_1\circ\xi)=
\da(i_2\circ i_1)\circ\da(\xi)\simeq\id$$ and
$$\da(\xi)\circ\da(i_2)=\da(i_1\circ\xi)\circ\da(i_2)=
\da(i_1)\circ\da(\xi\circ i_2)=\da(i_1)\circ\da(\xi)\simeq\id.$$
By a~symmetric argument the inclusion map
$j_2:\phi^{-1}(\ti\gamma) \hookrightarrow\phi^{-1}(\gamma')$
induces a~homotopy equivalence as well. Composing, we get that
$\da(i):\da(\phi^{-1}(\ti\gamma))\hookrightarrow\da(\phi^{-1}(\gamma))$
is a~homotopy equivalence.

In the special case $\gamma=(q>q')$ we get that
\[\da(\phi^{-1}(q))\simeq\da(\phi^{-1}(\gamma))\simeq
\da(\phi^{-1}(q')).\] Hence, since $Q$ is connected as a~poset, we get
$\da(\phi^{-1}(q))\simeq \da(\phi^{-1}(q'))$ for any $q,q'\in Q$.

Finally, the existence of the fibration homotopy long exact
sequence~(\ref{eq:piseq}) follows from (\ref{eq:bdgam}) and Quillen's
Theorem B, see \cite[p.\ 89]{Qu}. \qed

\begin{rem}
We shall not use Proposition~\ref{prop_ABC} (2) in this paper. We have
proved it here as a~result which is interesting on its own right and
might be useful for other computations. A more general version of
Proposition~\ref{prop_ABC} (1) was proved in~\cite{Ba}, 
see also~\cite{SZ}.
\end{rem}

\section{Complexes of homomorphisms from complete graphs}
\label{s4}

\subsection{The neighborhood complex and $\thom(K_2,G)$.} $\,$
\vspace{5pt}

We are now ready to define the neighborhood complex $\cn(G)$ and show
that it is homotopy equivalent to $\thom(K_2,G)$. The natural
advantage to working with the polyhedral complex $\thom(K_2,G)$
instead of the simplicial complex $\cn(G)$ is that $\thom(K_2,G)$
possesses a~natural free cellular $\zz$-action induced from the
swapping $\zz$-action on~$K_2$.

\begin{df}
  For an arbitrary graph $G$ the simplicial complex $\cn(G)$ is
  defined as follows: its vertices are all non-isolated vertices of
  $G$, and its simplices all the subsets of $V(G)$ which have a~common
  neighbor.
\end{df}
\nin In other words, the maximal simplices of $\cn(G)$ are $\tn(v)$,
for $v\in V(G)$.

\begin{prop} \label{propn}
$\thom(K_2,G)$ is homotopy equivalent to $\cn(G)$.
\end{prop}

\pr Let $P=\cp(\thom(K_2,G))$ and $Q=\cp(\cn(G))$. Consider
$\phi:P\rightarrow Q$ mapping the element $\eta:\{1,2\}\rightarrow
2^{V(G)}\setminus\emptyset$ to $\eta(1)\subseteq V(G)$. Clearly, the
vertices in $\eta(1)$ have all the vertices in $\eta(2)$ as their
neighbors, hence, since $\eta(2)\neq\emptyset$, $\phi$ is
well-defined. Let us show that $\phi$ induces homotopy equivalence
$\Delta(\phi):\Delta(P)\rightarrow\Delta(Q)$.

First, let $A\in Q$. We see that $\phi^{-1}(A)$ is the set of all
pairs $(A,B)$, $A,B\subseteq V(G)$, such that for all $x\in A$, and
$y\in B$, we have $(x,y)\in E(G)$. Clearly, $\phi^{-1}(A)$ has
a~maximal element $(A,\tn(A))$, so $\Delta(\phi^{-1}(A))$ is a~cone,
hence contractible.

Second, let us check the Condition $(B)$ of the
Proposition~\ref{prop_ABC}. Let $A\in Q$ and $(C,D)\in P$, such that
$\phi(C,D)=C\supseteq A$. Clearly $\tn(A)\supseteq\tn(C)\supseteq
D\neq\emptyset$. Then $\phi^{-1}(A)\cap P_{\leq(C,D)}=
\{(A,B)\,|\,B\subseteq D,B\neq\emptyset\}$. This poset has a~maximal
element $(A,D)$, since $D\subseteq\tn(A)$. In the notations of the
Proposition~\ref{prop_ABC} we have $(A,D)=\max((C,D),A)$.

Since Conditions $(A)$ and $(B)$ are satisfied, $\Delta(\phi)$ is
a~homotopy equivalence by Proposition~\ref{prop_ABC}. This shows that
$\bd(\thom(K_2,G))\simeq\bd(\cn(G))$, hence the result. \qed

\vspace{5pt} As Proposition \ref{propn} shows, the original complexes
$\cn(G)$ correspond to $K_2$-type obstructions to colorability. The
Lov\'asz' idea behind his Conjecture~\ref{loconj} was that the next
natural class of obstructions should come from the maps from odd
cycles $C_{2r+1}$ to our graph.

\subsection{$\thom(K_2,K_n)$ as a boundary complex of a polytope.} $\,$
\vspace{5pt}

Let $M_n$ denote the Minkowski sum
$$[-1/2,1/2]^n+[(-1/2,-1/2,\dots,-1/2),(1/2,1/2,\dots,1/2)],$$
where
$[-1/2,1/2]^n$ denotes the cube in ${\mathbb R}^n$ with vertices are
all points whose coordinates have the absolute value $1/2$. $M_n$ is
a~zonotope in ${\mathbb R}^n$. Its dual, $M_n^*$, is the polytope
associated to the hyperplane arrangement
$\ca=\{\ca_1,\dots,\ca_{n+1}\}$ defined by
$$\ca_i=\begin{cases}
(x_i=0),&\text{ for } 1\leq i\leq n;\\
(\sum_{j=1}^n x_j=0),&\text{ for } i=n+1.
\end{cases}$$

In the proof of the next proposition we identify each cell
$\eta:V(K_2)\ra 2^{V(K_n)}\sm\{\emptyset\}$ with the ordered pair
$(A,B)$ of non-empty subsets of $[n]$, by taking $A=\eta(1)$ and
$B=\eta(2)$.

\begin{figure}[hbt]
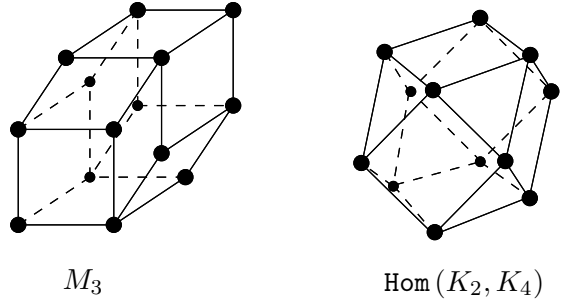

\begin{center}
  \begin{picture}(0,0)%
    \includegraphics{homk2kn.pstex}%
  \end{picture}%
  \input{homk2kn.pstex_t}%
  
\end{center}
\caption{Complex of 4-colorings of an edge and its dual.}
\label{fig:homk2kn}
\end{figure}

\begin{prop}\label{polyprop}
  $\thom(K_2,K_{n+1})$ is isomorphic as a~cell complex to the~boundary
  complex of $M_n^*$. The $\zz$-action on $\thom(K_2,K_{n+1})$,
  induced by the flip action of $\zz$ on $K_2$, corresponds under this
  isomorphism to the central symmetry.
\end{prop}

\pr Set $P=\cp(\thom(K_2,K_{n+1}))^{op}$. We shall see that $P$ is
isomorphic to the face poset of $M_n$, which we denote by~$Q={\mathcal
  F}(M_n)$. We shall denote the future isomorphism by~$\rho$.

First, note that faces of the cube $[-1/2,1/2]^n$ are encoded by
$n$-tuples of $1/2$, $-1/2$, and $*$, where $*$ denotes the coordinate
where the value can be chosen arbitrarily from the interval
$[-1/2,1/2]$. For an arbitrary $n$-tuple $x$, we let
$\supp(x)\subseteq[n]$ denote the set of the indices of coordinates
which are either non-zero, or are denoted with a~$*$. Additionally,
for an~arbitrary number $k$, we let $\supp(x,k)\subseteq[n]$ denote
the set of the indices of the coordinates which are equal to~$k$ (in
particular, they cannot be denoted with a~$*$).

Vertices of $M_n$ are labeled by all $n$-tuples of $1$, $-1$, and $0$,
such that $1$ and $-1$ are not present simultaneously, and not all the
coordinates are equal to~$0$, that is $v$ is a~vertex of $M_n$ if and
only if $v\in\{0,1\}^n$, or $v\in\{0,-1\}^n$, and $v\neq(0,\dots,0)$.
These vertices correspond to atoms in $P$ as follows:
\[
v\quad\stackrel{\rho}{\longleftrightarrow}\quad
\begin{cases}
(\supp(v),[n+1]\sm\supp(v)),&\text{ if  } v\in\{0,1\}^n; \\
([n+1]\sm\supp(v),\supp(v)),&\text{ if  } v\in\{0,-1\}^n.
\end{cases}
\]
Clearly, restricted to atoms, $\rho$ is a~bijection.

Those faces of $M_n$ which are contained in the closed star of
$(1,\dots,1)$ can be indexed by $f\in\{0,1,*\}^n$, where
$|\supp(f,1)|\geq 1$. Symmetrically, those faces of $M_n$ which are
contained in the closed star of $(-1,\dots,-1)$ can be indexed by
$f\in\{0,-1,*\}^n$, where $|\supp(f,-1)|\geq 1$. For these faces
$\rho$ can be defined as follows:
\[
f\quad\stackrel{\rho}{\longleftrightarrow}\quad
\begin{cases}
  (\supp(f,1),\supp(f,0)\cup\{n+1\}),&\text{ if  }
f\in\overline{\text{St}\,(1,\dots,1)}; \\
  (\supp(f,0)\cup\{n+1\},\supp(f,-1)),&\text{ if }
  f\in\overline{\text{St}\,(-1,\dots,-1)}.
\end{cases}
\]

Finally, we consider the faces of $M_n$ which are not in
$\overline{\text{St}\,(1,\dots,1)}\cup\overline{\text{St}\,(-1,\dots,-1)}$.
Each such face is a~convex hull of the union of two faces, $f\cup\ti
f$, such that $f\in\overline{\text{St}\,(1,\dots,1)}$, $\ti f\in
\overline{\text{St}\,(-1,\dots,-1)}$, with the condition that
$\supp(f,0)=\supp(\ti f,-1)$, $\supp(f,1)=\supp(\ti f,0)$. The element
of $P$ associated to such a~face under $\rho$ is
$(\supp(f,1),\supp(f,0))=(\supp(\ti f,0),\supp(\ti f,-1))$.

It is an easy exercise to check that $\rho$ defines a~poset
isomorphism between $P$ and~$Q$, which in turn induces the required
cell complex isomorphism.

Finally, a brief scanning through the definition of $\rho$ in
different cases reveals that $\rho$ is equivariant with respect to the
described $\zz$-actions on both sides. Hence the last part of the
proposition follows.
\qed
\vspace{5pt}

The cellular map $\phi$ defined in the Proposition~\ref{propn}, is in
this case going from the boundary of an~$(n-2)$-dimensional polytope
$M_n^*$ to the boundary of an~$(n-2)$-dimensional simplex. It would be
interesting to see whether it has interesting additional properties in
the context of zonotopes and also to find out what other graphs $G$
provide a connection to polytopes.

\begin{rem}
  $\thom(K_m,K_n)$ can be viewed as a~deleted product of simplices,
  see~\cite{Ma}. In this context it is well-known, probably due to
  van~Kampen, that for $m=2$ it is a boundary of a~polytope.
\end{rem}

\subsection{The homotopy type of $\thom(K_m,K_n)$.} $\,$ \label{ss4.3}
\vspace{5pt}

We can still get a~fairly detailed information about the topology of
the spaces of homomorphisms between complete graphs in general.

\begin{prop} \label{pr_chom}
  $\thom(K_m,K_n)$ is homotopy equivalent to a wedge of
  $(n-m)$-dimensional spheres.
\end{prop}
\pr We use induction on $m$ and on $n-m$. The base is provided by the
cases $\thom(K_1,K_n)$, which is a~simplex with $n$ vertices, hence
contractible, and $\thom(K_n,K_n)$, which consists of $n!$ points,
that is a~wedge of $n!-1$ spheres of dimension~$0$. We assume now that
$m\geq 2$ and $n\geq m+1$.

For $i\in[m]$ let $A_i$ be the subcomplex of $\thom(K_m,K_n)$ defined
by:
\[
A_i=\{\eta:[m]\ra 2^{[n]}\sm\{\emptyset\}\,|\, n\notin\eta(j), \text{
  for }j\in[m],j\neq i\}.
\]
Since any two vertices of $K_m$ are connected by an edge, $n$ cannot
be in $\eta(i_1)\cap\eta(i_2)$, for $i_1\neq i_2$. This implies that
$\bigcup_{i=1}^m A_i=\thom(K_m,K_n)$.

Clearly, for any $i\neq j$, $i,j\in[m]$, we have
\[
A_i\cap A_j=\{\eta:[m]\ra 2^{[n]}\sm\{\emptyset\}\,|\,
n\notin\eta(k),\text{ for all }k\in[m]\},
\]
so $A_i\cap A_j$ is isomorphic to $\thom(K_m,K_{n-1})$, hence, by
induction, it is $(n-m-2)$-connected.

We shall now see that each $A_i$ is $(n-m-1)$-connected. Since all
$A_i$'s are isomorphic to each other, it is enough to consider $A_1$.
Let us describe a~partial matching on $\cp(A_1)$. For
$\eta\in\cp(A_1)$, such that $n\notin\eta(1)$, we set
$\mu(\eta):=\tilde\eta$, defined by:
\[
\tilde\eta(i)=\begin{cases}
\eta(1)\cup\{n\},& \text{ for }i=1;\\
\eta(i),&\text{ for }i=2,3,\dots,m.
\end{cases}
\]
Obviously, this is an acyclic matching and the critical cells form
a~subcomplex $\wti A\subseteq A$ defined by: $\eta\in\wti A$ in and
only if $\eta(1)=\{n\}$. Thus $\wti A=\thom(K_{m-1},K_{n-1})$. Since,
by the Proposition~\ref{DMT} $\wti A$ is homotopy equivalent to $A_1$,
and $\wti A$ is $(n-m-1)$-connected by the induction assumption, we
conclude that $A_i$ is $(n-m-1)$-connected for any~$i$.

It follows from \cite[Theorem 10.6(ii)]{Bj} that $\thom(K_m,K_n)$ is
$(n-m-1)$-connected. Since dimension of $\thom(K_m,K_n)$ is $n-m$, it
follows from \cite[(9.19)]{Bj} that $\thom(K_m,K_n)$ is homotopy
equivalent to a~wedge of spheres. \qed\vspace{5pt}

One can use the construction in the proof of the
Proposition~\ref{pr_chom} to count the number of the spheres in the
wedge. Let us say that $\thom(K_m,K_n)$ is homotopy equivalent to
a~wedge of $f(m,n)$ spheres. Let $S(-,-)$ denote the Stirling numbers
of the second kind, and $SF_k(x)=\sum_{n\geq k}S(n,k)x^n$ denote the
generating function for these numbers. It is well-known that 
$$SF_k(x)=x^k/(1-x)(1-2x)\dots(1-kx).$$ 
For $m\geq 1$, let $F_m(x)=\sum_{n\geq 1}f(m,n)x^n$ be the generating
function for the number of the spheres. Clearly, $F_1(x)=0$, and
$F_2(x)=x^2/(1-x)$.

\begin{prop}
The numbers $f(m,n)$ satisfy the following recurrence relation
\begin{equation}\label{eqrecf}
f(m,n)=mf(m-1,n-1)+(m-1)f(m,n-1),
\end{equation}
for $n>m\geq 2$; with the boundary values $f(n,n)=n!-1$, $f(1,n)=0$
for $n\geq 1$, and $f(m,n)=0$ for $m>n$.

The generating function $F_m(x)$ is given by the equation:
\begin{equation}\label{eqngenff}
F_m(x)=(m!\cdot x \cdot SF_{m-1}(x)-x^m)/(1+x).
\end{equation}

As a~consequence, the following non-recursive formulae are valid:
\begin{equation}\label{eqnrecf2}
f(m,n)=(-1)^{m+n+1}+m!(-1)^{n}\sum_{k=m}^{n}(-1)^k S(k-1,m-1),
\end{equation}
and 
\begin{equation}\label{eqnrecf}
f(m,n)=\sum_{k=1}^{m-1}(-1)^{m+k+1}{m\choose k+1} k^n,
\end{equation}
for $n\geq m\geq 1$.
\end{prop}

\pr Let $\chi(m,n)$ denote the non-reduced Euler characteristics of
the complexes $\thom(K_m,K_n)$, and, for $i=1,\dots,m$, let $A_i$ be
as in the proof of the Proposition~\ref{pr_chom}. Since
$\thom(K_m,K_n)=\bigcup_{i=1}^m A_i$, $A_i\cap
A_j=\thom(K_m,K_{n-1})$, for all $i\neq j$, and
$A_i\simeq\thom(K_{m-1},K_{n-1})$, for $i\in[m]$, by simple
inclusion-exclusion counting we conclude that
\begin{equation}\label{eqrecchi}
\chi(m,n)=m\chi(m-1,n-1)-(m-1)\chi(m,n-1),
\end{equation}
for $n>m\geq 2$, additionally $\chi(n,n)=n!$, $\chi(1,n)=1$, for
$n\geq 1$. Since $\chi(m,n)=1+(-1)^{m-n}f(m,n)$, a~simple computation
shows the validity of~the relation~\eqref{eqrecf}.

For $m\geq 1$, let $G_m(x)=\sum_{n\geq 1}\chi(m,n)x^n$. Multiplying
each side of the equation \eqref{eqrecchi} by $x^n$ and summing over
all $n$ yields $G_m(x)=m\cdot x\cdot G_{m-1}(x)-(m-1)\cdot x\cdot
G_m(x)$, implying
$$G_m(x)=\dfrac{mx}{1+(m-1)x}G_{m-1}(x),$$ for $m\geq 1$, and hence,
since $G_0(x)=1/(1-x)$, we get
\begin{multline}
G_m(x)=\dfrac{m!\cdot x^m}{(1-x)(1+x)(1+2x)\dots(1+(m-1)x)}= \\
m!\cdot x\cdot (-1)^{m-1}\cdot SF_{m-1}(-x)/(1-x),
\end{multline}
for $m\geq 0$. By multiplying the identity
$f(m,n)=(-1)^{m+n}(\chi(m,n)-1)$ with $x^n$ and summing over all
$n\geq m$, we get
\begin{multline}
F_m(x)=(-1)^m G_m(-x)-x^m/(1+x)= \\ 
(-1)^m\cdot m!\cdot(-x)\cdot(-1)^{m-1}\cdot SF_{m-1}(x)/(1+x)-x^m/(1+x)=\\
(m!\cdot x\cdot SF_{m-1}(x)-x^m)/(1+x).
\end{multline}

\eqref{eqnrecf2} follows from comparing the coefficients
in~\eqref{eqngenff}.

To prove~\eqref{eqnrecf} we see that it fits the boundary values and
satisfies the recurrence relation~\eqref{eqrecf}.
Verifying~\eqref{eqrecf} is straightforward, as is
checking~\eqref{eqnrecf} for $m=1$ and $m=2$.
Finally,~\eqref{eqnrecf} is seen for $n=m\geq 2$ by expanding the
expression $(e^x-1)^n\cdot e^{-x}$ by the binomial theorem and
comparing the coefficient of $x^n$ on both sides of the expansion.
\qed \vspace{5pt}

In particular, we have $f(2,n)=1$, for $n\geq 2$, $f(3,n)=2^n-3$, for
$n\geq 3$, $f(4,n)=3^n-4\cdot 2^n+6$, for $n\geq 4$,
$f(5,n)=4^n-5\cdot 3^n+10\cdot 2^n-10$, for $n\geq 5$.

\vspace{5pt}

We are now ready to prove the result announced in the beginning of
this paper.

\vspace{5pt}

{\noindent{\bf Proof of the Theorem~\ref{cothm}. }} If the graph $G$
is $(k+m-1)$-colorable, then there exists a~homomorphism $\phi:G\ra
K_{k+m-1}$. It induces a $\zz$-equivariant map
$$\phi^{K_m}:\thom(K_m,G)\ra \thom(K_m,K_{k+m-1}).$$
By the
Proposition~\ref{pr_chom} the space $\thom(K_m,K_{k+m-1})$ is homotopy
equivalent to a~wedge of $(k-1)$-spheres, hence, by dimensional
reasons, $\varpi_1^{k}(\thom(K_m,K_{k+m-1}))=0$. Since the
Stiefel-Whitney classes are functorial, the existence of the map
$\phi^{K_m}$ implies that $\varpi_1^{k}(\thom(K_m,G))=0$, which is a
contradiction to the assumption of the theorem. \qed

\section{Complexes of homomorphisms from forests and their complements
  to complete graphs}

\subsection{The minor neighbor reduction and its consequences.} $\,$
\vspace{5pt}

The next proposition, coupled with Propositions~\ref{polyprop} and
\ref{pr_chom}, will be our workhorse for computing concrete examples.

\begin{prop} \label{prop5.1}
  If $G$ and $H$ are graphs and $u$ and $v$ are distinct vertices of
  $G$, such that $\tn(v)\subseteq\tn(u)$, then the inclusion
  $i:G-v\hookrightarrow G$, resp.\ the homomorphism $\phi:G\ra G-v$
  mapping $v$ to $u$ and fixing other vertices, induce homotopy
  equivalences $i_H:\thom(G,H)\ra\thom(G-v,H)$, resp.\ 
  $\phi_H:\thom(G-v,H)\ra\thom(G,H)$.
\end{prop}
\pr Let us apply the Proposition \ref{prop_ABC} (1) for the cellular
map $i_H:\thom(G,H)\ra\thom(G-v,H)$. Take $\eta\in\cp(\thom(G-v,H))$,
$\eta:V(G)\sm\{v\}\rightarrow 2^{V(H)}\sm\{\emptyset\}$. We have
$$\cp((i_H)^{-1}(\eta))=\{\tau\in\cp(\thom(G,H))\,|\,\tau(w)=\eta(w),
\text{ for } w\neq v,w\in V(G)\}.$$ An element in
$\cp((i_H)^{-1}(\eta))$ is determined by its value on $v$. Take
$\tau\in\cp((i_H)^{-1}(\eta))$ such that
$$\tau(v)=\bigcap_{y\in\tn(v)}\tn(\eta(y))\supseteq
\bigcap_{y\in\tn(u)}\tn(\eta(y))\supseteq\eta(u)\neq\emptyset.$$
Clearly, $\tau$ is the maximal element of $\cp((i_H)^{-1}(\eta))$,
hence $\Delta(\cp((i_H)^{-1}(\eta)))$ is contractible, so the
Condition $(A)$ is satisfied.

Let us now check the Condition $(B)$. Take $\tau\in\cp(\thom(G,H))$,
$\eta\in\cp(\thom(G-v,H))$, such that for any $x\in V(G)\sm\{v\}$ we
have $\tau(x)\supseteq\eta(x)$. The set
$\cp((i_H)^{-1}(\eta))\cap\cp(\thom(G,H))_{\leq\tau}$ consists of all
$\nu\in\cp(\thom(G,H))$, such that for any $x\in V(G)$ we have
$\tau(x)\supseteq\nu(x)$, and for any $x\in V(G)\sm\{v\}$ we have
$\eta(x)=\nu(x)$. Thus, it has a~maximal element defined by:
\[\nu(x)=\begin{cases} \eta(x),& \text{ for } x\neq v,\,\,x\in V(G);\\
\tau(x),& \text{ for } x=v.
\end{cases}\]

Conditions $(A)$ and $(B)$ being satisfied, we now get that
$\bd(i_H)$, hence also $i_H$, is a~homotopy equivalence.

To see that $\phi_H$ is also a~homotopy equivalence note first that
$i_H\circ\phi_H=\id_{\thom(G-v,H)}$. Let $j$ be the homotopy inverse
of $i_H$, then $\phi_H\circ i_H\simeq j\circ i_H\circ\phi_H\circ
i_H=j\circ i_H\simeq\id_{\thom(G,H)}$.
\qed \vspace{5pt}

If $G$ is a graph, and $u,v\in V(G)$, $u\neq v$, such that
$\tn(v)\subseteq\tn(u)$, then we say that $G$ {\it reduces} to $G-v$.
We shall also say that $u$ {\it~dominates} $v$, or that $v$ is
{\it~dominated} by $u$. If in addition $\tn(v)\neq\tn(u)$ we say that
$u$ {\it~strongly dominates} $v$. We call $u$ and $v$ {\it~equivalent}
if $\tn(v)=\tn(u)$. The strong domination defines a~partial order
$P(G)$ on the set of equivalence classes. We call a graph
{\it~irreducible} if it does not reduce to any subgraph.

We note a simple, but useful property of the vertex domination: if
$u,v\in S\subseteq V(G)$, $u\neq v$, and $u$ dominates $v$ in $G$,
then $u$ dominates $v$ in $G[S]$. If $u$ strongly dominates $v$ in
$G$, it is not true in general that $u$ strongly dominates $v$ in
$G[S]$.

As already the example of the tree shows, the minimal subgraph of
$G$ to which it reduces is not unique. However the following
weaker version of uniqueness is true.

\begin{prop}
  Let $G$ be a graph and $S,S'\subseteq V(G)$, such that $G$ reduces
  both to $G[S]$ and to $G[S']$, and both $G[S]$ and $G[S']$ are
  irreducible, then $G[S]$ is isomorphic to $G[S']$.
\end{prop}

\pr We prove the statement by induction on the number of vertices
in~$G$. If $|V(G)|=1$, then $S=S'=V(G)$, so the result is trivially
true. Assume now that $|V(G)|\geq 2$.

Choose $M\subseteq V(G)$ containing exactly one vertex from each
maximal equivalence class in $P(G)$, and no other vertices. If
$M=V(G)$, then $G$ is irreducible, so we can assume that $M\neq V(G)$.
Let us show that there exists $\wti S\subseteq M$, such that $G$
reduces to $G[\wti S]$, and $G[S]$ is isomorphic to $G[\wti S]$.

Assume that no such $\wti S$ exists. Consider all the reduction
sequences $(\ti v_1,\dots,\ti v_{|V(G)|-|S|})$ leading from $G$ to
a~graph isomorphic to $G[S]$. Set $\{\ti v_i\}_{i\in I}:=M\cap\{\ti
v_1,\dots,\ti v_{|V(G)|-|S|}\}$, and choose the sequence which
minimizes $\sum_{i\in I}(|V(G)|-i)$. Denote this sequence by
$(w_1,\dots,w_{|V(G)|-|S|})$.

Set $\wti S:=V(G)\sm\{w_1,\dots,w_{|V(G)|-|S|}\}$, and $\{w_i\}_{i\in
  I}:=M\cap\{w_1,\dots,w_{|V(G)|-|S|}\}$. If each vertex of $G[\wti
S]$ is either in $\wti S\cap M$ or is dominated in $G[\wti S]$ by some
vertex in $\wti S\cap M$, then, since $G[\wti S]$ is irreducible, we
conclude that $\wti S\subseteq M$, yielding a~contradiction.

Thus we may pick the {\it smallest} $i$, such that there exists
$v\in\wti S\sm M$, which is not dominated by any vertex of
$M\sm\{w_1,\dots,w_i\}$ in $G_i=G-\{w_1,\dots,w_i\}$. By the choice of
$M$, and what is said above, we have $i\in[|V(G)|-|S|]$. Clearly,
since $v$ was dominated by some vertex of $M\sm\{w_1,\dots,w_i\}$ in
$G_{i-1}=G-\{w_1,\dots,w_{i-1}\}$, we have that $w_i\in M$, and $w_i$
is the only vertex of $M\sm\{w_1,\dots,w_i\}$ which dominates $v$ in
$G_{i-1}$. In particular, $w_i$ itself is not dominated by any other
vertex of $M\sm\{w_1,\dots,w_i\}$ in $G_{i-1}$.

By the choice of $i$, every vertex in $G_{i-1}$, which is not in
$M\sm\{w_1,\dots,w_i\}$, is dominated by some vertex in
$M\sm\{w_1,\dots,w_i\}$, hence $w_i$ is not strongly dominated by any
other vertex. Since $G_{i-1}\ra G_{i-1}-\{w_i\}=G_i$ is a~legal
reduction, there must exist a~vertex $w$ equivalent to $w_i$ in
$G_{i-1}$. We have $w\notin M$, since either $w=v$, or $w$
dominates~$v$.

Consider a~graph isomorphism $\varphi:G_{i-1}\ra G_{i-1}$, which swaps
the vertices $w_i$ and $w$, and fixes every other vertex. It is easy
to see that $(w_1,\dots,w_{i-1},\varphi(w_i),\varphi(w_{i+1}),\dots,
\varphi(w_{|V(G)|-|S|}))$ is a~legal reduction sequence leading from
$G$ to $G[\widehat S]$, such that $G[\widehat S]$ is isomorphic to
$G[S]$.

Furthermore, since removal of $w_i\in M$ was either replaced by or
swapped with the removal of $w\notin M$, the invariant, which we
minimized over the sequences, is actually smaller for this sequence
than for $(w_1,\dots,w_{|V(G)|-|S|})$. This is again a~contradiction.

Finally, consider the case $S,S'\subseteq M$. Since $|M|<|V(G)|$, we
can use the induction assumption to prove the theorem, as long as we
can show that $G[M]$ reduces to $G[S]$ and to $G[S']$. By the argument
above, we can choose $S$ so that, if $(w_1,\dots,w_{|V(G)|-|S|})$ is
the reduction sequence leading to $G[S]$, and $\{w_i\}_{i\in I}=
M\cap\{w_1,\dots,w_{|V(G)|-|S|}\}$, then, for any
$i=1,\dots,|V(G)|-|S|$, every vertex in $V(G)\sm\{w_1,\dots,w_i\}$ is
dominated by some vertex from $M\sm\{w_1,\dots,w_i\}$ in
$G-\{w_1,\dots,w_i\}$. It is then immediate that
$\{w_{i_1},\dots,w_{i_t}\}$ is the reduction sequence from $G[M]$ to
$G[S]$, where $I=\{i_1,\dots,i_t\}$, $i_1<\dots<i_t$. 

Indeed, for any $i\in I$, $w_i$ is dominated by some vertex in
$G-\{w_1,\dots,w_{i-1}\}$, hence it is dominated by some vertex from
$M\sm\{w_1,\dots,w_{i-1}\}$ in $G-\{w_1,\dots,w_{i-1}\}$. It follows
that $w_i$ is dominated by some vertex in $G[M\sm\{w_j\,|\,j\in
I,j<i\}]$, allowing to reduce the latter graph to
$G[M\sm\{w_j\,|\,j\in I,j\leq i\}]$. \qed \vspace{5pt}

For future reference we explicitly state the following consequence of
the Proposition~\ref{prop5.1}.

\begin{crl}\label{equicrl}
  Let $G$ be a~graph, and $S\subseteq V(G)$, such that $G$ reduces to
  $G[S]$. Assume $S$ is $\Gamma$-invariant for some
  $\Gamma\subseteq\aut(G)$. Then the inclusion $i:G[S]\hra G$ induces
  a~$\Gamma$-invariant homotopy equivalence
  $i_H:\thom(G,H)\ra\thom(G[S],H)$ for an~arbitrary graph~$H$.
\end{crl}

Note also, that the Proposition~\ref{prop5.1} cannot be generalized to
encompass arbitrary graph homomorphisms $\phi$ of $G$ onto $H$, where
$H$ is a~subgraph of $G$, and $\phi$ is identity on $H$. As an example
in the subsection~\ref{ss2.3} showed
$\thom(C_6,K_3)\not\simeq\thom(K_2,K_3)$ despite of the existence of
the folding map of $C_6$ onto $K_2$.

\subsection{The homotopy type of $\thom(F,K_n)$ and $\thom(\ovr F,K_n)$.}
$\,$ \vspace{5pt}

Next, we use the Proposition~\ref{prop5.1} to compute homotopy types
of the complexes of maps from finite forests to complete graphs.

\begin{prop} \label{homtree}
  If $T$ is a~tree with at least one edge, then the map
  $i_{K_n}:\thom(T,K_n)\ra\thom(K_2,K_n)$ induced by any inclusion
  $i:K_2\hookrightarrow T$ is a~homotopy equivalence, in particular
  $\thom(T,K_n)\simeq S^{n-2}$. As a~consequence, if $F$ is a~forest,
  and $T_1,\dots,T_k$ are all its connected components consisting of
  at least 2 vertices, then $\thom(F,K_n)\simeq\prod_{i=1}^k S^{n-2}$.
\end{prop}
\pr Let $T$ be a~tree with $k$ vertices, $k\geq 2$. Note the~general
fact, that if $v$ is a~leaf of a~tree, $u$ is the~vertex adjacent to
$v$, and $w\neq v$ is a~vertex adjacent to $u$, then
$\tn(w)\supseteq\tn(v)=\{u\}$, hence $T$ reduces to $T-v$.

Let us now number the vertices $v_1,\dots,v_k$ so that for any
$i\in[k-1]$, $v_i$ is a~leaf in $T-\{v_{i+1},\dots,v_k\}$. By the
previous observation
\[
T\ra T-\{v_k\}\ra T-\{v_{k-1},v_k\}\ra\dots\ra T-\{v_3,\dots,v_k\}
=T[\{v_1,v_2\}]=K_2
\]
is a~valid reduction sequence. Thus the first part of the statement follows
by the Proposition~\ref{prop5.1}.

That $\thom(T,K_n)\simeq S^{n-2}$ follows from the
Proposition~\ref{polyprop}. Finally, the formula for the homotopy type
of $\thom(F,K_n)$ follows from (3) in the subsection~\ref{ss2.4}.
\qed\vspace{5pt}

Let $S^n_a$ denote the $n$-dimensional sphere equipped with the
antipodal action of $\zz$; in the same way $S^n_t$ denotes the
$n$-dimensional sphere equipped with the trivial action of $\zz$.

Given two spaces $X$ and $Y$ with $\zz$-action, we let $X\simeq_\zz Y$
denote the $\zz$-equivariant homotopy equivalence.

\begin{prop} \label{homtree2}
  Let $T$ be a~tree with at least one edge and a~$\zz$-action
  determined by an~invertible graph homomorphism $\gamma:T\ra T$. If
  $\gamma$ flips an~edge in $T$, then $\thom(T,K_n)\simeq_\zz
  S^{n-2}_a$, otherwise $\thom(T,K_n)\simeq_\zz S^{n-2}_t$.
\end{prop}
\pr Assume $\gamma$ flips an~edge, that is there exist $a,b\in V(G)$,
such that $(a,b)\in E(G)$, $\gamma(a)=b$, and $\gamma(b)=a$. By the
Corollary~\ref{equicrl} the inclusion map $i:T[\{a,b\}]\hra T$ induces
a~$\zz$-equivariant homotopy equivalence
$\thom(T,K_n)\simeq_\zz\thom(K_2,K_n)$, where the last space has the
natural $\zz$-action induced by the flipping $\zz$-action on $K_2$. By
the Proposition~\ref{polyprop} we get $\thom(T,K_n)\simeq_\zz
S^{n-2}_a$.

Assume now, there is no edge flipped by $\gamma$. Since $T$ is
a~contractible finite CW complex (the topology is generated by fixing
homeomorphisms between edges of $T$ and a~standard unit interval) it
follows from \cite[p.\ 257]{Bre} that $\gamma$ must have~a fixed
point. Denote this point by~$x$. Clearly, either $x$ is a~vertex of
$T$ or $x$ is the middle-point of some edge $e\in E[G]$. In the latter
case, if the edge is not fixed pointwise, then it is flipped, which
contradicts our assumptions on~$\gamma$.

Thus we found a~vertex $v\in V(G)$ fixed by $\gamma$. If there exists
$e=(a,b)\in E(G)$, such that $\gamma(a)=a$, $\gamma(b)=b$, then
$i:T[\{a,b\}]\hra T$ induces a~$\zz$-equivariant homotopy equivalence
$\thom(T,K_n)\simeq_\zz\thom(K_2,K_n)$, where the $\zz$-action on the
last space is the trivial one. It follows that $\thom(T,K_n)\simeq_\zz
S^{n-2}_t$.

Finally consider the case when there is no edge in $T$ which is fixed
by $\gamma$ pointwise. Let $u$ be any vertex of $T$ adjacent to $v$,
and let $w=\gamma(u)\neq u$. Since the set $\{u,w,v\}$ is
$\gamma$-invariant, we see by the Corollary~\ref{equicrl} that the
inclusion map $i:T[\{u,w,v\}]\hra T$ induces a~$\zz$-equivariant
homotopy equivalence
$i_H:\thom(T,K_n)\ra\thom(T[\{u,w,v\}],K_n)=\thom(L_3,K_n)$, where the
$\zz$-action on the last space is induced from the $\zz$-action on
$L_3$ which swaps $u$ and~$w$.

Let $\phi:L_3\ra K_2$ be any of the two
elements of $\chom(L_3,K_2)$. Clearly, $\phi$ is $\zz$-equivariant
with the $\zz$-action on $K_2$ being trivial. This shows that
$\phi_H:\thom(L_3,K_n)\ra\thom(K_2,K_n)\simeq_\zz S^{n-2}_t$ is
a~$\zz$-equivariant homotopy equivalence, which finishes the proof.
\qed
\vspace{5pt}

Since taking the unlooped complement reverses neighbor set inclusions,
we see that $G$ reduces if and only if $\ovr{G}$ reduces. The next
proposition describes what happens if $G$ is a~forest.

\begin{prop} \label{homcompf}
  If $F$ is a~forest, then $\thom(\ovr{F},K_n)\simeq\thom(K_m,K_n)$,
  where $m$ is the maximal cardinality of an~independent set in~$F$.
\end{prop}
\pr We use induction on the number of edges in $F$. If
$E(F)=\emptyset$, then $\ovr{F}=K_{|V(F)|}$, the maximal cardinality
of an~independent set in~$F$ is $|V(F)|$, and the statement is
obvious. So assume $|E(F)|\geq 1$.

Let $v\in V(F)$ be an arbitrary leaf, and let $u\in V(F)$ be the
vertex adjacent to $v$. We have $\tn_{\ovr{F}}(u)\subseteq
V(F)\sm\{u,v\}=\tn_{\ovr{F}}(v)$. Hence $\ovr{F}$ reduces to
$\ovr{F}-u$. Clearly, $\ovr{F}-u=\ovr{F-u}$, so by combining the
induction assumption with the Proposition~\ref{prop5.1} we get
$\thom(\ovr{F},K_n)\simeq\thom(K_{\wti m},K_n)$, where $\wti m$ is the
maximal cardinality of an~independent set in~$F-u$.

Let $I$ be an independent set in $F$ of maximal cardinality. Either $u$
or $v$ must be in $I$, since otherwise $I\cup\{v\}$ is independent,
and larger than $I$. If $u\in I$, then $(I\sm\{u\})\cup\{v\}$ is also
an independent set in $F$ of maximal cardinality. Either way, we have
an independent set $J$ in $F$ of maximal cardinality containing~$v$.
Since any independent set in $F-u$ is also independent in $F$, and $J$
is independent in $F-u$, we can conclude that $m=\wti m$, hence the
result.
\qed
\vspace{5pt}

It follows from the Proposition~\ref{pr_chom} that
$\thom(\ovr{F},K_n)$ is homotopy equivalent to a~wedge of
$(n-m)$-dimensional spheres.

\section{Notes added in proof.}

After the appearance of the initial series of papers, \cite{BK1,BK2},
and this paper, there has been a~lot of further research activities
pertaining to $\thom$-complexes. We refer the reader to the recent
survey~\cite{Koz5}.

We update here on a few concrete developments which have taken place
by the time of the publication of this article.
\begin{itemize}
\item The Conjecture \ref{conj1} has been proved, see \cite{CK2}. \vskip5pt

\item The Proposition \ref{propn} has been strengthened. It has been 
proved that $\thom(K_2,G)$ and $\cn(G)$ have the same simple homotopy
type. See \cite{Koz4} for the proof, and a~description of an~explicit
formal deformation. \vskip5pt

\item The Proposition \ref{prop5.1} has been strengthened. It has been 
proved that the folds are admissible on both sides of the
$\thom(-,-)$, and that one can replace the homotopy equivalences by
formal deformations, see \cite{Koz3} for an~explicit construction.
\end{itemize}

\end{document}